\newif\ifJOURNAL
\JOURNALfalse
\newif\ifarXiv
\arXivfalse
\newif\ifWP
\WPfalse
\newif\ifFULL
\FULLfalse
\newif\ifSemiFULL
\SemiFULLfalse
\newif\ifLATIN
\LATINfalse

\arXivtrue

\ifFULL\SemiFULLtrue\fi

\ifarXiv\LATINtrue\fi	

\newif\ifnotJOURNAL	
\notJOURNALtrue
\ifJOURNAL\notJOURNALfalse\fi

\newif\ifnotarXiv	
\notarXivtrue
\ifarXiv\notarXivfalse\fi

\newif\ifTR		
\TRfalse
\ifarXiv\TRtrue\fi
\ifWP\TRtrue\fi
\newif\ifnotTR
\notTRtrue
\ifarXiv\notTRfalse\fi
\ifWP\notTRfalse\fi

\newif\ifnotFULL	
\notFULLtrue
\ifFULL\notFULLfalse\fi

\newif\ifnotLATIN	
\notLATINtrue
\ifLATIN\notLATINfalse\fi

\ifJOURNAL
  \newcommand{\CTI}{vovk:arXiv0712.1275}
  \newcommand{\CTII}{vovk:arXiv0712.1483}
\fi
\ifarXiv
  \newcommand{\CTI}{vovk:arXiv0712.1275}
  \newcommand{\CTII}{vovk:arXiv0712.1483}
\fi
\ifWP
  \newcommand{\CTI}{GTP24}
  \newcommand{\CTII}{GTP25}
\fi

\ifnotLATIN
  \newcommand{\Takeuchi}{takeuchi:2004}
\fi
\ifLATIN
  \newcommand{\Takeuchi}{takeuchi:2004latin}
\fi

\ifJOURNAL
\documentclass{article}
\usepackage{amsmath,amsthm,amsfonts,amssymb,latexsym}
\newcommand{\Extra}[1]{}
\fi

\ifarXiv
\documentclass{article}
\usepackage{amsmath,amsthm,amsfonts,amssymb,latexsym}
\newcommand{\Extra}[1]{}
\fi

\ifWP
\documentclass[toc]{gtarticle}
\usepackage{amsmath,amsthm,amsfonts,amssymb,latexsym,epsfig}
\renewcommand{\Extra}[1]{}
\fi

\ifFULL
\usepackage{color}
\renewcommand{\Extra}[1]{\blue{#1}}

\newcommand{\blue}[1]{\textcolor{blue}{#1}}
\newcommand{\bluebegin}{\begingroup\color{blue}}
\newcommand{\blueend}{\endgroup}

\fi

\allowdisplaybreaks[4]

\newcommand{\Vladimir}{Vladimir}
\newcommand{\DOT}{.}

\ifnotLATIN
  \input{OT2enc.def}
  
  \usepackage{CJK}
\fi

\newcommand{\st}{\mathrel{|}}		

\newcommand{\dd}{d}			

\newcommand{\K}{\mathcal{K}}		

\newcommand{\AAA}{\mathcal{A}}		
\newcommand{\FFF}{\mathcal{F}}		
\newcommand{\Normal}{\mathcal{N}}	

\DeclareMathOperator{\III}{\mathbb{I}}		

\newcommand{\bbbp}{\mathbb{P}}		

\DeclareMathOperator{\UpProb}{\overline{\bbbp}}		
\DeclareMathOperator{\LowProb}{\underline{\bbbp}}	

\DeclareMathOperator{\UpProbMod}{\smash{\overline{\bbbp}}\vphantom{\bbbp}^{\prime}}
\DeclareMathOperator{\LowProbMod}{\smash{\underline{\bbbp}}\vphantom{\bbbp}^{\prime}}

\newcommand{\bbbe}{\mathbb{E}}		
\DeclareMathOperator{\Expect}{\bbbe}
\DeclareMathOperator{\UpExpect}{\overline{\bbbe}}	
\DeclareMathOperator{\LowExpect}{\underline{\bbbe}}	

\DeclareMathOperator{\UpExpectMod}{\smash{\overline{\bbbe}}\vphantom{\bbbe}^{\prime}}
\DeclareMathOperator{\LowExpectMod}{\smash{\underline{\bbbe}}\vphantom{\bbbe}^{\prime}}

\newcommand{\bbbr}{\mathbb{R}}		

\theoremstyle{plain}
\newtheorem{theorem}{Theorem}
\newtheorem{proposition}{Proposition}
\newtheorem{corollary}{Corollary}
\newtheorem{lemma}{Lemma}

\theoremstyle{definition}
\newtheorem*{remark}{Remark}

\title{Game-theoretic Brownian motion}

\ifarXiv
\author{Vladimir Vovk\\
\texttt{vovk{\rm@}cs.rhul.ac.uk}\\
\texttt{http://vovk.net}}
\fi

\ifWP
\author{Vladimir Vovk}

\twodatestrue

\fi

\begin{document}

\maketitle

\begin{abstract}
  This paper suggests a perfect-information game,
  along the lines of L\'evy's characterization of Brownian motion,
  that formalizes the process of Brownian motion
  in game-theoretic probability.
  This is perhaps the simplest situation
  where probability emerges in a non-stochastic environment.
\end{abstract}

\section{Introduction}

This paper is part of the recent revival of interest
(see, e.g., \cite{dawid/vovk:1999,shafer/vovk:2001,\Takeuchi,%
kumon/etal:2007,horikoshi/takemura:2008,kumon/takemura:2008})
in game-theoretic probability.
It further develops the game-theoretic approach to continuous-time processes
along the lines of the papers \cite{takeuchi/etal:2007,\CTI,\CTII}.
Unlike those papers,
which only demonstrate the emergence of randomness-type properties
in a continuous trading protocol,
this paper gives an example where a full-fledged probability measure emerges
(although in a significantly more restrictive protocol).
Only the very simple case of ``game-theoretic Brownian motion'' is considered,
but we can expect that probabilities will also emerge
in less restrictive protocols.

The words such as ``positive'', ``negative'', ``before'', and ``after''
will be understood in the wide sense of $\ge$ or $\le$,
as appropriate;
when necessary, we will add the qualifier ``strictly''.

\ifarXiv
  The latest version of this working paper can be downloaded from the web site
  \texttt{http://probabilityandfinance.com}.
\fi

\section{Upper and lower probability}

We consider a perfect-information game between two players,
Reality and Sceptic.
Reality chooses a continuous function $\omega:[0,\infty)\to\bbbr$
with $\omega(0)=0$,
but before she announces her choice
Sceptic chooses his strategy of trading in two securities,
one with the price process $\omega(t)$, $t\in[0,\infty)$,
and the other with the price process $\omega^2(t)-t$, $t\in[0,\infty)$.
This game, which we call the \emph{L\'evy game}
(after \cite{levy:1965}, Theorem 18.6),
is formalized as follows.

Let $\Omega$ be the set of all continuous functions
$\omega:[0,\infty)\to\bbbr$ with $\omega(0)=0$.
For each $t\in[0,\infty)$,
$\FFF_t$ is defined to be the $\sigma$-algebra
generated by the functions
$\omega\in\Omega\mapsto\omega(s)$,
$s\in[0,t]$,
and $\FFF_{\infty}:=\vee_t\FFF_t$;
we will often write $\FFF$ for $\FFF_{\infty}$.
A \emph{process} $S$ is a family of functions
$S_t:\Omega\to[-\infty,\infty]$, $t\in[0,\infty)$,
each $S_t$ being $\FFF_t$-measurable;
we only consider processes
with lower continuous (often continuous) \emph{sample paths}
$t\mapsto S_t(\omega)$.
An \emph{event} is an element of the $\sigma$-algebra $\FFF$.
Stopping times $\tau:\Omega\to[0,\infty]$
w.r.\ to the filtration $(\FFF_t)$
and the corresponding $\sigma$-algebras $\FFF_{\tau}$
are defined as usual.
We simplify $\omega(\tau(\omega))$ and $S_{\tau(\omega)}(\omega)$
to $\omega(\tau)$ and $S_{\tau}(\omega)$,
respectively;
the argument $\omega$ will often be omitted
in other cases as well.

The class of allowed strategies for Sceptic is defined in two steps.
An \emph{elementary betting strategy} $G$
consists of an increasing sequence of stopping times
$\tau_1\le\tau_2\le\cdots$
and, for each $n=1,2,\ldots$,
a pair of bounded $\FFF_{\tau_{n}}$-measurable functions,
$M_n$ and $V_n$.
It is required that, for any $\omega\in\Omega$,
$\lim_{n\to\infty}\tau_n(\omega)=\infty$.
To such $G$ and an \emph{initial capital} $c\in\bbbr$
corresponds the \emph{elementary capital process}
\begin{multline}\label{eq:elementary-capital}
  \K^{G,c}_t(\omega)
  :=
  c
  +
  \sum_{n=1}^{\infty}
  \biggl(
    M_n(\omega)
    \bigl(
      \omega(\tau_{n+1}\wedge t)-\omega(\tau_n\wedge t)
    \bigr)\\
    +
    V_n(\omega)
    \Bigl(
      \bigl(
        \omega^2(\tau_{n+1}\wedge t)
        -
        (\tau_{n+1}\wedge t)
      \bigr)
      -
      \bigl(
        \omega^2(\tau_n\wedge t)
        -
        (\tau_n\wedge t)
      \bigr)
    \Bigr)
  \biggr),\\
  t\in[0,\infty)
\end{multline}
(with the zero terms in the sum ignored).
The numbers $M_n(\omega)$ and $V_n(\omega)$ will be called Sceptic's \emph{stakes}
(on $\omega(t)$ and $\omega^2(t)-t$, respectively)
chosen at time $\tau_n$,
and $\K^{G,c}_t(\omega)$ will sometimes be referred to
as Sceptic's capital at time $t$;
we may also say that Sceptic \emph{bets} $M_n(\omega)$ on $\omega(t)$
and $V_n(\omega)$ on $\omega^2(t)-t$ at time $\tau_n$.
(We are following standard probability textbooks,
such as \cite{williams:1991}, Chapter 10,
in using gambling rather than financial terminology.)

A \emph{positive capital process} is any process $S$
that can be represented in the form
\begin{equation}\label{eq:positive-capital}
  S_t(\omega)
  :=
  \sum_{n=1}^{\infty}
  \K^{G_n,c_n}_t(\omega),
\end{equation}
where the elementary capital processes $\K^{G_n,c_n}_t(\omega)$
are required to be positive, for all $t$ and $\omega$,
and the positive series $\sum_{n=1}^{\infty}c_n$ is required to converge.
The sum (\ref{eq:positive-capital}) is always positive
but allowed to take value $\infty$.
Since $\K^{G_n,c_n}_0(\omega)=c_n$ does not depend on $\omega$,
$S_0(\omega)$ also does not depend on $\omega$
and will often be abbreviated to $S_0$.

The \emph{upper probability} of a set $E\subseteq\Omega$
is defined as
\begin{equation}\label{eq:upper-probability}
  \UpProb(E)
  :=
  \inf
  \bigl\{
    S_0
    \bigm|
    \forall\omega\in\Omega:
    \liminf_{t\to\infty}
    S_t(\omega)
    \ge
    \III_E(\omega)
  \bigr\},
\end{equation}
where $S$ ranges over the positive capital processes
and $\III_E$ stands for the indicator function of $E$.
(The $\liminf_{t\to\infty}$ can be replaced by $\sup_{t\in[0,\infty)}$
in this definition:
see \cite{\CTI}, Lemma 1.)
The \emph{lower probability} of $E\subseteq\Omega$ is
\begin{equation*}
  \LowProb(E)
  :=
  1-\UpProb(E^c),
\end{equation*}
where $E^c:=\Omega\setminus E$ stands for the complement of $E$.

\begin{remark}
  Our definition of a positive capital process
  corresponds to the intuitive picture
  where Sceptic divides his initial capital
  into a sequence of independent accounts,
  with a prudent (not risking bankruptcy) elementary betting strategy
  applied to each account.
  On the other hand,
  we could make the definition of upper probability
  more similar to the standard definition of expectation
  for positive random variables:
  a positive capital process could be equivalently defined
  as the limit of an increasing sequence of positive elementary capital processes
  with uniformly bounded initial capitals.
\end{remark}

\section{Emergence of the Wiener measure}

It is obvious that upper probability is countably (in particular, finitely) subadditive:
\begin{lemma}\label{lem:subadditivity}
  For any sequence of subsets $E_1,E_2,\ldots$ of $\Omega$,
  \begin{equation*}
    \UpProb
    \left(
      \bigcup_{n=1}^{\infty}
      E_n
    \right)
    \le
    \sum_{n=1}^{\infty}
    \UpProb(E_n).
  \end{equation*}
\end{lemma}
\noindent
Therefore,
$\UpProb$ is an outer measure in Carath\'eodory's sense.
Recall that a set $A\subseteq\Omega$ is \emph{$\UpProb$-measurable}
if, for each $E\subseteq\Omega$,
\begin{equation}\label{eq:measurable}
  \UpProb(E)
  =
  \UpProb(E\cap A)
  +
  \UpProb(E\cap A^c).
\end{equation}
A standard result
(see \cite{caratheodory:1914}, Sections 1--11,
or, e.g., \cite{kallenberg:2002}, Theorem 2.1)
shows that the family $\AAA$ of all $\UpProb$-measurable sets
forms a $\sigma$-algebra
and that the restriction of $\UpProb$ to $\AAA$ is a probability measure on $(\Omega,\AAA)$.

\begin{theorem}\label{thm:emergence}
  Each event $A\in\FFF$ is $\UpProb$-measurable,
  and the restriction of $\UpProb$ to $\FFF$
  coincides with the Wiener measure $W$ on $(\Omega,\FFF)$.
  In particular,
  $\UpProb(A)=\LowProb(A)=W(A)$ for each $A\in\FFF$.
\end{theorem}
\noindent
The rest of this paper is devoted to proving this result.

\section{Statement in terms of expectation}
\label{sec:expectation}

A \emph{capital process} is a process of the form $S-C$,
where $S$ is a positive capital process and $C$ is a real constant.
The \emph{upper expectation} of a bounded functional $F:\Omega\to\bbbr$ is
\begin{equation*}
  \UpExpect(F)
  :=
  \inf
  \bigl\{
    S_0
    \bigm|
    \forall\omega\in\Omega:
    \liminf_{t\to\infty}
    S_t(\omega)
    \ge
    F(\omega)
  \bigr\},
\end{equation*}
where $S$ ranges over the capital processes.
(This generalizes upper probability:
$\UpProb(E)=\UpExpect(\III_E)$ for all $E\subseteq\Omega$.)
Theorem \ref{thm:emergence} will immediately follow from the following result:
\begin{theorem}\label{thm:expectation}
  If $F$ is a bounded $\FFF$-measurable functional on $\Omega$,
  \begin{equation}\label{eq:equality}
    \UpExpect(F)
    =
    \int_{\Omega} F(\omega) W(\dd\omega).
  \end{equation}
\end{theorem}

Indeed, let us deduce Theorem \ref{thm:emergence} from Theorem \ref{thm:expectation}.
Let $A\in\FFF$.
Since the inequality $\le$ in (\ref{eq:measurable})
follows from Lemma \ref{lem:subadditivity},
to show that $A\in\AAA$
we are only required to show
\begin{equation}\label{eq:to-show}
  \UpProb(E\cap A)
  +
  \UpProb(E\cap A^c)
  \le
  \UpProb(E)
  +
  \epsilon
\end{equation}
for each $\epsilon>0$.
Fix such an $\epsilon$.

Let $S$ be a positive capital process such that 
$S_0<\UpProb(E)+\epsilon$ and
$
  \forall\omega\in\Omega:
  \liminf_{t\to\infty}
  S_t(\omega)
  \ge
  \III_E(\omega)
$.
Set
$
  F(\omega)
  :=
  1
  \wedge
  \liminf_{t\to\infty}
  S_t(\omega)
$,
so that $F$ is a bounded 
(taking values in $[0,1]$)
$\FFF$-measurable functional satisfying
$\UpExpect(F)<\UpProb(E)+\epsilon$,
$\III_{E\cap A}\le F\III_A$,
and $\III_{E\cap A^c}\le F\III_{A^c}$
(the last two inequalities follow from $\III_E\le F$).
Therefore,
it suffices to prove
\begin{equation*}
  \UpExpect(F\III_A)
  +
  \UpExpect(F\III_{A^c})
  \le
  \UpExpect(F).
\end{equation*}
This immediately follows from
$
  F\III_A
  +
  F\III_{A^c}
  =
  F
$
and Theorem \ref{thm:expectation}.

The remaining statements of Theorem~\ref{thm:emergence}
are obvious corollaries of Theorem~\ref{thm:expectation}:
for $A\in\FFF$,
\begin{align*}
  \UpProb(A)
  &=
  \UpExpect(\III_A)
  =
  \int_{\Omega}
  \III_A
  \dd W
  =
  W(A),\\
  \LowProb(A)
  &=
  1 - \UpProb(A^c)
  =
  1 - W(A^c)
  =
  W(A).
\end{align*}

\section{Coherence and its application}

In this and the following two sections we will prove
some auxiliary results
that will be needed in the proof of Theorem~\ref{thm:expectation}.
This section's results, however,
also have considerable substantive significance:
they show that our definitions are ``free of contradiction''.
(In fact, these definitions have been chosen to make this easy.)

The following result says that the L\'evy game is \emph{coherent},
in the sense that $\UpProb(\Omega)=1$
(i.e.,
no positive capital process increases its value between time $0$ and $\infty$
by more than a positive constant for all $\omega\in\Omega$).
\begin{proposition}\label{prop:coherence}
  $\UpProb(\Omega)=1$.
\end{proposition}
\begin{proof}
  If $\omega$ is generated as sample path of (measure-theoretic) Brownian motion,
  any positive elementary capital process will be a positive continuous local martingale
  (since, by the optional sampling theorem,
  every partial sum in (\ref{eq:elementary-capital})
  will be a continuous martingale),
  and so it suffices to apply the maximal inequality for positive supermartingales
  to the partial sums corresponding to a given positive capital process.
  \ifSemiFULL

    In the rest of this proof I will check carefully
    the martingale property.
    I will only consider the sum
    \begin{equation*}
      \sum_{n=1}^{\infty}
      M_n(\omega)
      \bigl(
        \omega(\tau_{n+1}\wedge t)-\omega(\tau_n\wedge t)
      \bigr)
    \end{equation*}
    assuming that $\omega$ is a continuous martingale
    (e.g., $\omega^2(t)-t$ in the old notation).
    We can rewrite this sum as
    \begin{equation*}
      \sum_{n=1}^{\infty}
      h_n(\omega)
      \bigl(
        \omega(t)-\omega(\tau_n\wedge t)
      \bigr)
    \end{equation*}
    for some bounded $\FFF_{\tau_n}$-measurable $h_n$
    (namely, $h_n:=M_n-M_{n-1}$, with $M_0:=0$).
    Our goal is to check that each addend
    \begin{equation}\label{eq:addend}
      h_n(\omega)
      \bigl(
        \omega(t)-\omega(\tau_n\wedge t)
      \bigr)
    \end{equation}
    is indeed a martingale.
    (See \cite{\CTI} for a marginally more complicated argument
    establishing directly that each addend in (\ref{eq:elementary-capital})
    is a martingale.)
    For each $t\in[0,\infty)$,
    (\ref{eq:addend}) is integrable
    by the boundedness of $h_n$ and the optional sampling theorem
    (see, e.g., \cite{revuz/yor:1999}, Theorem II.3.2).
    We only need to prove, for $0<s<t$, that
    (omitting, until the end of the proof, the argument $\omega$ and ``a.s.'')
    \begin{equation}\label{eq:to-prove}
      \Expect
      \left(
        h_n
        \bigl(
          \omega(t)-\omega(\tau_n\wedge t)
        \bigr)
        \mid
        \FFF_s
      \right)
      =
      h_n
      \bigl(
        \omega(s)-\omega(\tau_n\wedge s)
      \bigr).
    \end{equation}
    We will check this equality on two $\FFF_s$-measurable events separately:
    \begin{description}
    \item[$\{\tau_n\le s\}$:]
      We need to check
      \begin{equation*}
        \Expect
        \left(
          h_n
          \bigl(
            \omega(t)-\omega(\tau_n)
          \bigr)
          \III_{\{\tau_n\le s\}}
          \mid
          \FFF_s
        \right)
        =
        h_n
        \bigl(
          \omega(s)-\omega(\tau_n)
        \bigr)
        \III_{\{\tau_n\le s\}}.
      \end{equation*}
      Since $h_n\III_{\{\tau_n\le s\}}$ is bounded and $\FFF_s$-measurable
      (its $\FFF_s$-measurability follows, e.g.,
      from Lemma 1.2.15 in \cite{karatzas/shreve:1991} and the monotone-class theorem),
      it suffices to check
      \begin{equation*}
        \Expect
        \left(
          \bigl(
            \omega(t)-\omega(\tau_n)
          \bigr)
          \III_{\{\tau_n\le s\}}
          \mid
          \FFF_s
        \right)\\
        =
        \bigl(
          \omega(s)-\omega(\tau_n)
        \bigr)
        \III_{\{\tau_n\le s\}}.
      \end{equation*}
      Since $\omega(\tau_n)\III_{\{\tau_n\le s\}}$ is $\FFF_s$-measurable,
      it suffices to check
      \begin{equation*}
        \Expect
        \left(
          \omega(t)
          \III_{\{\tau_n\le s\}}
          \mid
          \FFF_s
        \right)
        =
        \omega(s)
        \III_{\{\tau_n\le s\}}.
      \end{equation*}
      The stronger equality
      \begin{equation*}
        \Expect
        \left(
          \omega(t)
          \mid
          \FFF_s
        \right)
        =
        \omega(s)
      \end{equation*}
      is part of the definition of a martingale.
    \item[$\{s<\tau_n\}$:]
      We are required to prove
      \begin{equation*}
        \Expect
        \left(
          h_n
          \bigl(
            \omega(t)-\omega(\tau_n\wedge t)
          \bigr)
          \III_{\{s<\tau_n\}}
          \mid
          \FFF_s
        \right)
        =
        0,
      \end{equation*}
      but we will prove more:
      \begin{equation*}
        \Expect
        \left(
          h_n
          \bigl(
            \omega(t)-\omega(\tau_n\wedge t)
          \bigr)
          \III_{\{s<\tau_n\}}
          \mid
          \FFF_{s\vee\tau_n\wedge t}
        \right)
        =
        0
      \end{equation*}
      ($s\vee x\wedge t$ being a shorthand
      for $(s\vee x)\wedge t$ or, equivalently, $s\vee(x\wedge t)$).
      Since the event $\{\tau_n\le t\}$,
      being equivalent to $\tau_n\le s\vee\tau_n\wedge t$,
      is $\FFF_{s\vee\tau_n\wedge t}$-measurable
      (see \cite{karatzas/shreve:1991}, Lemma~1.2.16),
      it is sufficient to prove
      \begin{equation}\label{eq:remains}
        \Expect
        \left(
          h_n
          \bigl(
            \omega(t)-\omega(\tau_n\wedge t)
          \bigr)
          \III_{\{s<\tau_n\le t\}}
          \mid
          \FFF_{s\vee\tau_n\wedge t}
        \right)
        =
        0
      \end{equation}
      and
      \begin{equation*}
        \Expect
        \left(
          h_n
          \bigl(
            \omega(t)-\omega(\tau_n\wedge t)
          \bigr)
          \III_{\{t<\tau_n\}}
          \mid
          \FFF_{s\vee\tau_n\wedge t}
        \right)
        =
        0.
      \end{equation*}
      The second equality is obvious,
      so our task has reduced to proving the first, (\ref{eq:remains}).
      Since $h_n\III_{\{\tau_n\le t\}}=h_n\III_{\{\tau_n\le s\vee\tau_n\wedge t\}}$
      is bounded and $\FFF_{s\vee\tau_n\wedge t}$-measurable,
      (\ref{eq:remains}) reduces to
      \begin{equation*}
        \Expect
        \left(
          \bigl(
            \omega(t)-\omega(\tau_n\wedge t)
          \bigr)
          \III_{\{s<\tau_n\le t\}}
          \mid
          \FFF_{s\vee\tau_n\wedge t}
        \right)
        =
        0,
      \end{equation*}
      which is the same thing as
      \begin{equation*}
        \Expect
        \left(
          \bigl(
            \omega(t)-\omega(s\vee\tau_n\wedge t)
          \bigr)
          \III_{\{s<\tau_n\le t\}}
          \mid
          \FFF_{s\vee\tau_n\wedge t}
        \right)
        =
        0.
      \end{equation*}
      The last equality follows from the $\FFF_{s\vee\tau_n\wedge t}$-measurability
      of the event $\{s<\tau_n\le t\}=\{s<s\vee\tau_n\wedge t\le t\}$
      (see \cite{karatzas/shreve:1991}, Lemma 1.2.16)
      and the special case
      \begin{multline*}
        \Expect
        \left(
          \bigl(
            \omega(t)-\omega(s\vee\tau_n\wedge t)
          \bigr)
          \mid
          \FFF_{s\vee\tau_n\wedge t}
        \right)\\
        =
        \omega(s\vee\tau_n\wedge t) - \omega(s\vee\tau_n\wedge t)
        =
        0.
      \end{multline*}
      of the optional sampling theorem.
      \qedhere
    \end{description}
  \fi
\end{proof}

The \emph{lower expectation} of a bounded functional $F:\Omega\to\bbbr$
is defined as
\begin{equation*}
  \LowExpect(F)
  :=
  -\UpExpect(-F).
\end{equation*}
\begin{corollary}\label{cor:expectation}
  For every bounded functional on $\Omega$,
  $\LowExpect(F)\le\UpExpect(F)$.
\end{corollary}
\begin{proof}
  Suppose $\LowExpect(F)>\UpExpect(F)$ for some $F$;
  by the definition of $\LowExpect$,
  this would mean that $\UpExpect(F)+\UpExpect(-F)<0$.
  Since $\UpExpect$ is finitely subadditive,
  this would imply $\UpExpect(0)<0$,
  which is equivalent to $\UpProb(\Omega)<1$
  and, therefore, impossible by Proposition \ref{prop:coherence}.
\end{proof}

Specializing Corollary \ref{cor:expectation}
to the indicator functions of subsets of $\Omega$,
we obtain:
\begin{corollary}\label{cor:probability}
  For every set $A\subseteq\Omega$,
  $\LowProb(A)\le\UpProb(A)$.
\end{corollary}

\noindent
In the special case where $A$ is $\UpProb$-measurable
(see (\ref{eq:measurable})),
Proposition \ref{prop:coherence} implies $\UpProb(A)=\LowProb(A)$.

The facts that we have just established
allow us to simplify our goal:
(\ref{eq:equality}) will follow from
\begin{equation}\label{eq:inequality}
  \UpExpect(F)
  \le
  \int_{\Omega} F(\omega) W(\dd\omega).
\end{equation}
Indeed,
(\ref{eq:inequality}) implies
\begin{equation*}
  \LowExpect(F)
  =
  -\UpExpect(-F)
  \ge
  -\int_{\Omega} (-F(\omega)) W(\dd\omega)
  =
  \int_{\Omega} F(\omega) W(\dd\omega),
\end{equation*}
and so Corollary \ref{cor:expectation} implies
\begin{equation*}
  \UpExpect(F)
  =
  \LowExpect(F)
  =
  \int_{\Omega} F(\omega) W(\dd\omega).
\end{equation*}

\section{Modified L\'evy game}

Theorem \ref{thm:expectation} (as well as Theorem \ref{thm:emergence}) also holds,
and will be slightly easier to prove,
for the natural modification of the L\'evy game
in which Sceptic is allowed to bet
on $\omega(t)-\omega(\tau)$ and $(\omega(t)-\omega(\tau))^2-(t-\tau)$
at any time $\tau$.
Formally,
we obtain the definition of upper probability
in the modified L\'evy game
when we replace (\ref{eq:elementary-capital}) by
\begin{multline*}
  \K^{G,c}_t(\omega)
  :=
  c
  +
  \sum_{n=1}^{\infty}
  \biggl(
    M_n(\omega)
    \bigl(
      \omega(\tau_{n+1}\wedge t)-\omega(\tau_n\wedge t)
    \bigr)\\
    +
    V_n(\omega)
    \left(
      \bigl(
        \omega(\tau_{n+1}\wedge t)-\omega(\tau_n\wedge t)
      \bigr)^2
      -
      \bigl(
        (\tau_{n+1}\wedge t)-(\tau_n\wedge t)
      \bigr)
    \right)
  \biggr);
\end{multline*}
we will say that the corresponding elementary betting strategy
bets $M_n(\omega)$
(or stakes $M_n(\omega)$ units) on $\omega(t)-\omega(\tau_n)$
and bets $V_n(\omega)$ (or stakes $V_n(\omega)$ units)
on $(\omega(t)-\omega(\tau_n))^2-(t-\tau_n)$
at time $\tau_n$.
The rest of the definition is as before:
positive capital processes are defined by (\ref{eq:positive-capital})
and upper probability is then defined by (\ref{eq:upper-probability}).
The definitions of lower probability and upper and lower expectation
also carry over to the modified L\'evy game.

The L\'evy game and modified L\'evy game are very close,
as can be seen from the identity
\begin{multline}\label{eq:identity}
  M_n(\omega)
  \bigl(
    \omega(\tau_{n+1}\wedge t)-\omega(\tau_n\wedge t)
  \bigr)\\
  +
  V_n(\omega)
  \left(
    \bigl(
      \omega(\tau_{n+1}\wedge t)-\omega(\tau_n\wedge t)
    \bigr)^2
    -
    \bigl(
      (\tau_{n+1}\wedge t)-(\tau_n\wedge t)
    \bigr)
  \right)\\
  =
  \bigl(
    M_n(\omega)
    -
    2\omega(\tau_n\wedge t)
    V_n(\omega)
  \bigr)
  \bigl(
    \omega(\tau_{n+1}\wedge t)-\omega(\tau_n\wedge t)
  \bigr)\\
  +
  V_n(\omega)
  \Bigl(
    \bigl(
      \omega^2(\tau_{n+1}\wedge t)
      -
      (\tau_{n+1}\wedge t)
    \bigr)
    -
    \bigl(
      \omega^2(\tau_n\wedge t)
      -
      (\tau_n\wedge t)
    \bigr)
  \Bigr).
\end{multline}
However, the absence of an upper bound on $\omega$
prevents us from asserting
that the two games lead to the same notion of upper probability.
(Remember that the stakes are required to be bounded,
which is implicitly used in the proof of Proposition~\ref{prop:coherence}.)
If there is a risk of confusion,
we will write $\UpProbMod$, $\LowProbMod$, $\UpExpectMod$, and $\LowExpectMod$
instead of $\UpProb$, $\LowProb$, $\UpExpect$, and $\LowExpect$, respectively,
for the modified L\'evy game.

It is obvious that Lemma~\ref{lem:subadditivity}
continues to hold for the modified L\'evy game.
This is also true about Proposition~\ref{prop:coherence}
and Corollaries~\ref{cor:expectation}--\ref{cor:probability}
(although this fact is not used in this paper
outside this section).
As already mentioned,
Theorem \ref{thm:emergence} and Theorem \ref{thm:expectation}
also hold for the modified L\'evy game:
we will prove (\ref{eq:inequality}) directly for this game,
and the reduction to (\ref{eq:inequality}) depended on arguments
of general nature,
not involving the specifics of the L\'evy game.

\section{Tightness in the modified L\'evy game}

The set $\Omega$ is equipped with the standard metric
\begin{equation}\label{eq:metric}
  \rho(\omega_1,\omega_2)
  :=
  \sum_{n=1}^{\infty}
  2^{-n}
  \sup_{t\in[0,n]}
  \left(
    \lvert \omega_1(t)-\omega_2(t) \rvert
    \wedge
    1
  \right),
\end{equation}
which makes $\Omega$ a complete separable metric space
with Borel $\sigma$-algebra $\FFF$.
In this topology,
$\LowProbMod$ is tight:
\begin{lemma}\label{lem:tight}
  For each $\alpha>0$ there exists a compact set $K\subseteq\Omega$
  such that $\LowProbMod(K)\ge1-\alpha$.
\end{lemma}

For $\omega\in\Omega$ and $T\in(0,\infty)$,
the \emph{modulus of continuity of $\omega$ on $[0,T]$}
is defined as
\begin{equation*}
  m^T_{\delta}(\omega)
  :=
  \sup_{s,t\in[0,T]:\lvert s-t\rvert\le\delta}
  \lvert\omega(s)-\omega(t)\rvert,
  \quad
  \delta>0.
\end{equation*}
The proof of Lemma \ref{lem:tight} will be based on the following result.
\begin{lemma}\label{lem:modulus}
  For each $\alpha>0$ and $T>0$,
  \begin{equation}\label{eq:modulus}
    \LowProbMod
    \left\{
      \forall\delta>0:
      m^T_{\delta}
      \le
      157\,
      \alpha^{-1/2}
      T^{3/8}
      \delta^{1/8}
    \right\}
    \ge
    1-\alpha.
  \end{equation}
\end{lemma}
\noindent
Of course,
Theorem \ref{thm:emergence} and its counterpart for the modified L\'evy game
will imply much subtler results
than Lemma~\ref{lem:modulus} and its counterpart for the L\'evy game,
such as L\'evy's modulus of continuity formula.
It is interesting, however,
that this section's results
do not require that Sceptic should be allowed to bet on $\omega(t)-\omega(\tau)$
(i.e., he can achieve his goal even if he is required
to always choose $M_n(\omega):=0$).

\begin{proof}[Proof of Lemma \ref{lem:modulus}.]
  We will use the method of \cite{vovk:1993forecasting},
  pp.~213--216.
  For each $n=1,2,\ldots$,
  divide the time interval $[0,T]$ into $2^n$ equal subintervals of length $2^{-n}T$.
  Fix, for a moment, an $n$,
  and set
  $
    \beta=\beta_n
    :=
    \left(2^{1/4}-1\right)
    2^{-n/4}
    \alpha
  $
  (where $2^{1/4}-1$ is the normalizing constant
  ensuring that the $\beta_n$ sum to $\alpha$)
  and
  \begin{equation*}
    \omega_i
    :=
    \omega(i2^{-n}T),
    \quad
    i=0,1,\ldots,2^n.
  \end{equation*}
  With lower probability at least $1-\beta/2$,
  \begin{equation}\label{eq:2-variation}
    \sum_{i=1}^{2^n}
    (\omega_i - \omega_{i-1})^2
    \le
    2T/\beta
  \end{equation}
  (since there is a positive elementary capital process taking value
  $
    T
    +
    \sum_{i=1}^{j}
    (\omega_i - \omega_{i-1})^2
    -
    j2^{-n}T
  $
  at time $j2^{-n}T$,
  $j=0,1,\ldots,2^n$,
  and this elementary capital process
  will make $2T/\beta$ at time $T$ out of initial capital $T$
  if (\ref{eq:2-variation}) fails to happen).
  For each $\omega\in\Omega$,
  define
  \begin{equation*}
    J(\omega)
    :=
    \left\{
      i=1,\ldots,2^n:
      \lvert\omega_i-\omega_{i-1}\rvert\ge\epsilon
    \right\},
  \end{equation*}
  where $\epsilon=\epsilon(\delta,n)$ will be chosen later.
  It is clear that $\lvert J(\omega)\rvert\le 2T/\beta\epsilon^2$
  on the set (\ref{eq:2-variation}).
  Consider the elementary betting strategy
  that bets $1$ on $(\omega(t)-\omega(\tau))^2-(t-\tau)$
  at each time $\tau\in[(i-1)2^{-n}T,i2^{-n}T]$ with $i\in J(\omega)$
  when $\lvert\omega(\tau)-\omega_{i-1}\rvert=\epsilon$
  for the first time during $[(i-1)2^{-n}T,i2^{-n}T]$
  and gets rid of the stake at time $i2^{-n}T$.
  This strategy will make at least $\epsilon^2$ out of $(2T/\beta\epsilon^2)2^{-n}T$
  provided both event (\ref{eq:2-variation})
  and event
  \begin{equation*}
    \exists i\in\{1,\ldots,2^n\}:
    \lvert\omega_i-\omega_{i-1}\rvert\ge2\epsilon
  \end{equation*}
  happen.
  (And we can make the corresponding elementary capital process positive
  by allowing Sceptic to bet at most $2T/\beta\epsilon^2$ times.)
  This corresponds to making at least $1$ out of $(2T/\beta\epsilon^4)2^{-n}T$.
  Solving the equation $(2T/\beta\epsilon^4)2^{-n}T=\beta/2$
  gives $\epsilon=(4T^22^{-n}/\beta^2)^{1/4}$.
  Therefore,
  \begin{multline}\label{eq:max}
    \max_{i=1,\ldots,2^n}
    \lvert\omega_i-\omega_{i-1}\rvert
    \le
    2\epsilon
    =
    2(4T^22^{-n}/\beta^2)^{1/4}\\
    =
    2^{3/2}
    \left(
      2^{1/4} - 1
    \right)^{-1/2}
    \alpha^{-1/2}
    T^{1/2}
    2^{-n/8}
  \end{multline}
  with lower probability at least $1-\beta$.
  By the countable subadditivity of upper probability
  (Lemma~\ref{lem:subadditivity}),
  (\ref{eq:max}) holds for all $n=1,2,\ldots$
  with lower probability at least $1-\sum_n\beta_n=1-\alpha$.

  Intervals of the form $[(i-1)2^{-n}T,i2^{-n}T]$,
  for $n\in\{1,2,\ldots\}$ and $i\in\{1,2,3,\ldots,2^n\}$,
  will be called \emph{dyadic}.
  Given an interval $[s,t]$ of length at most $\delta>0$ in $[0,T]$,
  we can cover its interior
  (without covering any points in its complement)
  by adjacent dyadic intervals with disjoint interiors
  such that, for some $m\in\{1,2,\ldots\}$:
  there are between one and two dyadic intervals of length $2^{-m}T$;
  for $i=m+1,m+2,\ldots$,
  there are at most two dyadic intervals of length $2^{-i}T$
  (start from finding the point in $[s,t]$ of the form $2^{-k}T$
  with the smallest possible $k$ and cover $(s,2^{-k}T]$ and $[2^{-k}T,t)$
  by dyadic intervals in the greedy manner).
  Combining (\ref{eq:max}) and
  $
    2^{-m}T\le\delta
  $,
  we obtain:
  \begin{align*}
    m^T_{\delta}(\omega)
    &\le
    2
    \times
    2^{3/2}
    \left(
      2^{1/4} - 1
    \right)^{-1/2}
    \alpha^{-1/2}
    T^{1/2}\\
    &\quad\times
    \left(
      2^{-m/8}
      +
      2^{-(m+1)/8}
      +
      2^{-(m+2)/8}
      +
      \cdots
    \right)\\
    &=
    2^{5/2}
    \left(
      2^{1/4} - 1
    \right)^{-1/2}
    \left(
      1 - 2^{-1/8}
    \right)^{-1}
    \alpha^{-1/2}
    T^{1/2}
    2^{-m/8}\\
    &\le
    2^{5/2}
    \left(
      2^{1/4} - 1
    \right)^{-1/2}
    \left(
      1 - 2^{-1/8}
    \right)^{-1}
    \alpha^{-1/2}
    T^{1/2}
    (\delta/T)^{1/8}\\
    &=
    2^{5/2}
    \left(
      2^{1/4} - 1
    \right)^{-1/2}
    \left(
      1 - 2^{-1/8}
    \right)^{-1}
    \alpha^{-1/2}
    T^{3/8}
    \delta^{1/8},
  \end{align*}
  which is stronger than (\ref{eq:modulus}).
\end{proof}

Now we can prove the following elaboration of Lemma \ref{lem:tight},
which will also be used in the next section.
\begin{lemma}\label{lem:super-modulus}
  For each $\alpha>0$,
  \begin{equation}\label{eq:super-modulus}
    \LowProbMod
    \left\{
      \forall T\ge1\;
      \forall\delta>0:
      m^T_{\delta}
      \le
      560\,
      \alpha^{-1/2}
      T^{1/2}
      \delta^{1/8}
    \right\}
    \ge
    1-\alpha.
  \end{equation}
\end{lemma}
\begin{proof}
  Replacing $\alpha$ in (\ref{eq:modulus}) by
  $\alpha_T:=(1-2^{-1/4})T^{-1/4}\alpha$ for $T=1,2,4,8,\ldots$
  (where $1-2^{-1/4}$ is the normalizing constant
  ensuring that the $\alpha_T$ sum to $\alpha$ over $T$),
  we obtain
  \begin{equation*}
    \LowProbMod
    \left\{
      \forall\delta>0:
      m^T_{\delta}
      \le
      157\,
      (1-2^{-1/4})^{-1/2}
      \alpha^{-1/2}
      T^{1/2}
      \delta^{1/8}
    \right\}
    \ge
    1-(1-2^{-1/4})T^{-1/4}\alpha.
  \end{equation*}
  The countable subadditivity of upper probability now gives
  \begin{multline*}
    \LowProbMod
    \left\{
      \forall T\in\{1,2,4,\ldots\}\;
      \forall\delta>0:
      m^T_{\delta}
      \le
      157\,
      (1-2^{-1/4})^{-1/2}
      \alpha^{-1/2}
      T^{1/2}
      \delta^{1/8}
    \right\}\\
    \ge
    1-\alpha,
  \end{multline*}
  which in turn gives
  \begin{equation*}
    \LowProbMod
    \left\{
      \forall T\ge1\;
      \forall\delta>0:
      m^T_{\delta}
      \le
      157\,
      (1-2^{-1/4})^{-1/2}
      \alpha^{-1/2}
      (2T)^{1/2}
      \delta^{1/8}
    \right\}
    \ge
    1-\alpha,
  \end{equation*}
  which is stronger than (\ref{eq:super-modulus}).
\end{proof}

Lemma \ref{lem:tight} immediately follows
from Lemma \ref{lem:super-modulus}
and the Arzel\`a--Ascoli theorem
(as stated in \cite{karatzas/shreve:1991}, Theorem 2.4.9).

Inequality (\ref{eq:2-variation}) will also be useful in the next section;
the following lemma packages it in a convenient form.
\begin{lemma}\label{lem:super-2-variation}
  For each $\alpha>0$,
  \begin{multline}\label{eq:super-2-variation}
    \LowProbMod
    \Biggl\{
      \forall T\in\{1,2,4,\ldots\}\;
      \forall n\in\{1,2,\ldots\}:
      \sum_{i=1}^{2^n}
      \Bigl(
        \omega(i2^{-n}T)
        -
        \omega((i-1)2^{-n}T)
      \Bigr)^2\\
      \le
      46\,
      \alpha^{-1}
      T^2
      2^{n/16}
    \Biggr\}
    \ge
    1-\alpha.
  \end{multline}
\end{lemma}
\begin{proof}
  Replacing $\beta/2$ in (\ref{eq:2-variation})
  with $2^{-1}(2^{1/16}-1)T^{-1}2^{-n/16}\alpha$,
  where $T$ ranges over $\{1,2,4,\ldots\}$
  and $n$ over $\{1,2,\ldots\}$,
  we obtain
  \begin{multline*}
    \LowProbMod
    \Biggl\{
      \sum_{i=1}^{2^n}
      \Bigl(
        \omega(i2^{-n}T)
        -
        \omega((i-1)2^{-n}T)
      \Bigr)^2\\
      \le
      2(2^{1/16}-1)^{-1}
      \alpha^{-1}
      T^2
      2^{n/16}
    \Biggr\}
    \ge
    1-2^{-1}(2^{1/16}-1)T^{-1}2^{-n/16}\alpha;
  \end{multline*}
  by the countable subadditivity of upper probability this implies
  \begin{multline*}
    \LowProbMod
    \Biggl\{
      \forall T\in\{1,2,4,\ldots\}\;
      \forall n\in\{1,2,\ldots\}:
      \sum_{i=1}^{2^n}
      \Bigl(
        \omega(i2^{-n}T)
        -
        \omega((i-1)2^{-n}T)
      \Bigr)^2\\
      \le
      2(2^{1/16}-1)^{-1}
      \alpha^{-1}
      T^2
      2^{n/16}
    \Biggr\}
    \ge
    1-\alpha,
  \end{multline*}
  which is stronger than (\ref{eq:super-2-variation}).
\end{proof}

\section{Proof of (\ref{eq:inequality}) for the modified L\'evy game}
\label{sec:proof}

To establish (\ref{eq:inequality})
(with $\UpExpect$ replaced by $\UpExpectMod$)
we only need to establish $\UpExpectMod(F)<\int F\dd W+\epsilon$
for a positive constant $\epsilon$.
We start from a series of reductions:
\begin{enumerate}
\item
  We can assume that $F$ is lower semicontinuous on $\Omega$.
  Indeed, if it is not,
  by the Vitali--Carath\'eodory theorem
  (see, e.g., \cite{rudin:1987short}, Theorem 2.24)
  for any compact $K\subseteq\Omega$
  there exists a lower semicontinuous function $G$ on $K$
  such that $G\ge F$ on $K$
  and $\int_K G \dd W \le \int_K F \dd W + \epsilon$.
  Without loss of generality we assume $\sup G\le\sup F$,
  and we extend $G$ to all of $\Omega$
  by setting $G:=\sup F$ outside $K$.
  Choosing $K$ with large enough $W(K)$
  (which can be done since the probability measure $W$ is tight:
  see, e.g., \cite{billingsley:1968}, Theorem 1.4),
  we will have $G\ge F$
  and $\int G \dd W \le \int F \dd W + 2\epsilon$.
  Achieving $S_0\le\int G\dd W+\epsilon$
  and $\liminf_{t\to\infty}S_t(\omega)\ge G(\omega)$,
  where $S$ is a capital process,
  will automatically achieve $S_0\le\int F\dd W + 3\epsilon$
  and $\liminf_{t\to\infty}S_t(\omega)\ge F(\omega)$.
\item
  We can further assume that $F$ is continuous on $\Omega$.
  Indeed, since each lower semicontinuous function on a metric space
  is a limit of an increasing sequence of continuous functions
  (see, e.g., \cite{engelking:1989}, Problem 1.7.15(c)),
  given a lower semicontinuous function $F$ on $\Omega$
  we can find a series of positive continuous functions $G^n$ on $\Omega$,
  $n=1,2,\ldots$,
  such that $\inf F+\sum_{n=1}^{\infty}G^n=F$.
  The sum $S$ of $\inf F$ and positive capital processes $S^1,S^2,\ldots$
  achieving $S^n_0\le\int G^n\dd W+2^{-n}\epsilon$
  and $\liminf_{t\to\infty}S^n_t(\omega)\ge G^n(\omega)$,
  $n=1,2,\ldots$,
  will achieve $S_0\le\int F\dd W + \epsilon$
  and $\liminf_{t\to\infty}S_t(\omega)\ge F(\omega)$.
\item
  We can further assume that $F$ depends on $\omega\in\Omega$
  only via $\omega|_{[0,T]}$ for some $T\in(0,\infty)$.
  Indeed, let us fix $\epsilon>0$ and prove
  $\UpExpectMod(F)\le\int F\dd W + C\epsilon$
  for some positive constant $C$
  assuming $\UpExpectMod(G)\le\int G\dd W$
  for all $G$ that depend on $\omega$
  only via $\omega|_{[0,T]}$ for some $T\in(0,\infty)$.
  Choose a compact set $K\subseteq\Omega$
  with $W(K)>1-\epsilon$ and $\LowProbMod(K)>1-\epsilon$
  (cf.\ Lemma~\ref{lem:tight}).
  Set
  $
    F^T(\omega)
    :=
    F(\omega^T)
  $,
  where $\omega^T$ is defined by $\omega^T(t):=\omega(t\wedge T)$
  and $T$ is sufficiently large in the following sense.
  Since $F$ is uniformly continuous on $K$
  and the metric is defined by (\ref{eq:metric}),
  $F$ and $F^T$ can be made arbitrarily close in $C(K)$
  (spaces $C(\ldots)$ are always equipped with the $\sup$ norm
  in this paper);
  in particular, let $\|F-F^T\|_{C(K)}<\epsilon$.
  Choose capital processes $S^0$ and $S^1$ such that
  $S^0_0\le\int F^T\dd W+\epsilon$,
  $\liminf_{t\to\infty}S^0_t(\omega)\ge F^T(\omega)$,
  $S^1_0\le\epsilon$,
  $\liminf_{t\to\infty}S^1_t(\omega)\ge\III_{K^c}(\omega)$.
  The sum $S:=S^0+2\sup\lvert F\rvert S^1+\epsilon$ will satisfy
  \begin{align*}
    S_0
    &\le
    \int F^T\dd W
    +
    (2\sup\lvert F\rvert+2)\epsilon
    \le
    \int_K F^T\dd W
    +
    (3\sup\lvert F\rvert+2)\epsilon\\
    &\le
    \int_K F\dd W
    +
    (3\sup\lvert F\rvert+3)\epsilon
    \le
    \int F\dd W
    +
    (4\sup\lvert F\rvert+3)\epsilon
  \end{align*}
  and
  \begin{equation*}
    \liminf_{t\to\infty}S_t(\omega)
    \ge
    F^T(\omega)
    +
    2\sup\lvert F\rvert
    \III_{K^c}(\omega)
    +
    \epsilon
    \ge
    F(\omega).
  \end{equation*}
  Without loss of generality,
  we assume $T\in\{1,2,4,\ldots\}$.
\item
  We can further assume that $F$ depends on $\omega$
  only via the values $\omega(iT/N)$,
  $i=1,\ldots,N$ (remember that $\omega(0)=0$),
  for some $N\in\{1,2,\ldots\}$.
  Indeed, let us fix $\epsilon>0$ and prove
  $\UpExpectMod(F)\le\int F\dd W + C\epsilon$
  for some positive constant $C$
  assuming $\UpExpectMod(G)\le\int G\dd W$
  for all $G$ that depend on $\omega$
  only via $\omega(iT/N)$, $i=1,\ldots,N$, for some $N$.
  Let $K\subseteq\Omega$ be the compact set in $\Omega$
  defined as
  $
    K
    :=
    \left\{
      \omega
      \st
      \forall\delta>0:
      m^T_{\delta}
      \le
      f(\delta)
    \right\}
  $
  for some $f:(0,\infty)\to(0,\infty)$ satisfying
  $\lim_{\delta\to0}f(\delta)=0$
  (cf.\ the Arzel\`a--Ascoli theorem)
  and chosen in such a way that
  $W(K)>1-\epsilon$ and $\LowProbMod(K)>1-\epsilon$.
  Let $g$ be a modulus of continuity of $F$ on $K$;
  we know that $\lim_{\delta\to0}g(\delta)=0$.
  Set
  $
    F_N(\omega)
    :=
    F(\omega_N)
  $,
  where $\omega_N$ is the piecewise linear function
  whose graph is obtained
  by joining the points $(iT/N,\omega(iT/N))$, $i=0,1,\ldots,N$,
  and $(\infty,\omega(T))$,
  and $N$ is so large that $g(f(T/N))\le\epsilon$.
  Since
  \begin{equation*}
    \omega\in K
    \enspace\Longrightarrow\enspace
    \left\|\omega-\omega_N\right\|_{C([0,T])}
    \le
    f(T/N)
    \enspace\Longrightarrow\enspace
    \rho(\omega,\omega_N)
    \le
    f(T/N)
  \end{equation*}
  (we assume, without loss of generality,
  that the graph of $\omega$ is horizontal over $[T,\infty)$),
  we have $\|F-F_N\|_{C(K)}\le\epsilon$.
  Choose capital processes $S^0$ and $S^1$ such that
  $S^0_0\le\int F_N\dd W+\epsilon$,
  $\liminf_{t\to\infty}S^0_t(\omega)\ge F_N(\omega)$,
  $S^1_0\le\epsilon$,
  $\liminf_{t\to\infty}S^1_t(\omega)\ge\III_{K^c}(\omega)$.
  The sum $S:=S^0+2\sup\lvert F\rvert S^1+\epsilon$ will satisfy
  \begin{align*}
    S_0
    &\le
    \int F_N\dd W
    +
    (2\sup\lvert F\rvert+2)\epsilon
    \le
    \int_K F_N\dd W
    +
    (3\sup\lvert F\rvert+2)\epsilon\\
    &\le
    \int_K F\dd W
    +
    (3\sup\lvert F\rvert+3)\epsilon
    \le
    \int F\dd W
    +
    (4\sup\lvert F\rvert+3)\epsilon
  \end{align*}
  and
  \begin{equation*}
    \liminf_{t\to\infty}S_t(\omega)
    \ge
    F_N(\omega)
    +
    2\sup\lvert F\rvert
    \III_{K^c}(\omega)
    +
    \epsilon
    \ge
    F(\omega).
  \end{equation*}
\item
  We can further assume that
  \begin{equation}\label{eq:F}
    F(\omega)
    =
    U
    \left(
      \omega(T/N),
      \omega(2T/N),
      \ldots,
      \omega(T)
    \right)
  \end{equation}
  where the function $U:\bbbr^{N}\to\bbbr$
  is not only continuous but also has compact support.
  (We will sometimes say that $U$ is the \emph{generator} of $F$.)
  Indeed, let us fix $\epsilon>0$ and prove
  $\UpExpectMod(F)\le\int F\dd W + C\epsilon$
  for some positive constant $C$
  assuming $\UpExpectMod(G)\le\int G\dd W$
  for all $G$ whose generator has compact support.
  Let $B_R$ be the open ball of radius $R$ and centred at the origin
  in the space $\bbbr^{N}$ with the $\ell_{\infty}$ norm.
  We can rewrite (\ref{eq:F})
  as $F(\omega)=U(s(\omega))$
  where $s:\Omega\to\bbbr^N$ reduces each $\omega\in\Omega$
  to
  $
    s(\omega)
    :=
    \left(
      \omega(T/N),
      \omega(2T/N),
      \ldots,
      \omega(T)
    \right)
  $.
  Choose $R$ so large that
  $W(s^{-1}(B_R))>1-\epsilon$ and $\LowProbMod(s^{-1}(B_R))>1-\epsilon$
  (the existence of such $R$ follows
  from the Arzel\`a--Ascoli theorem and Lemma~\ref{lem:tight}).
  Alongside $F$, whose generator is denoted $U$,
  we will also consider $F^*$ with generator
  \begin{equation*}
    U^*(\sigma)
    :=
    \begin{cases}
      U(\sigma) & \text{if $\sigma\in\overline{B_R}$}\\
      0 & \text{if $\sigma\in B^c_{2R}$};
    \end{cases}
  \end{equation*}
  in the remaining region $B_{2R}\setminus\overline{B_R}$,
  $U^*$ is defined arbitrarily
  (but making sure that $U^*$ is continuous
  and takes values in $[\inf U,\sup U]$;
  this can be done by the Tietze--Urysohn theorem,
  \cite{engelking:1989}, Theorem 2.1.8).
  Choose capital processes $S^0$ and $S^1$ such that
  $S^0_0\le\int F^*\dd W+\epsilon$,
  $\liminf_{t\to\infty}S^0_t(\omega)\ge F^*(\omega)$,
  $S^1_0\le\epsilon$,
  $\liminf_{t\to\infty}S^1_t(\omega)\ge\III_{s^{-1}(B_R^c)}(\omega)$.
  The sum $S:=S^0+2\sup\lvert F\rvert S^1$ will satisfy
  \begin{align*}
    S_0
    &\le
    \int F^*\dd W
    +
    (2\sup\lvert F\rvert+1)\epsilon
    \le
    \int_{s^{-1}(B_R)} F^*\dd W
    +
    (3\sup\lvert F\rvert+1)\epsilon\\
    &=
    \int_{s^{-1}(B_R)} F\dd W
    +
    (3\sup\lvert F\rvert+1)\epsilon
    \le
    \int F\dd W
    +
    (4\sup\lvert F\rvert+1)\epsilon
  \end{align*}
  and
  \begin{equation*}
    \liminf_{t\to\infty}S_t(\omega)
    \ge
    F^*(\omega)
    +
    2\sup\lvert F\rvert
    \III_{s^{-1}(B_R^c)}(\omega)
    \ge
    F(\omega).
  \end{equation*}
\item
  Since every continuous $U:\bbbr^{N}\to\bbbr$
  with compact support
  can be arbitrarily well approximated in $C(\bbbr^{N})$
  by an infinitely differentiable function with compact support
  (see, e.g., \cite{adams/fournier:2003short}, Theorem 2.29(d)),
  we can further assume that the generator $U$ of $F$
  is an infinitely differentiable function with compact support.
\item
  By Lemma~\ref{lem:tight},
  it suffices to prove that,
  given $\epsilon>0$ and a compact set $K$ in $\Omega$,
  some capital process $S\ge\inf F-1$
  with $S_0\le\int F \dd W + \epsilon$ achieves
  $\liminf_{t\to\infty}S_t(\omega)\ge F(\omega)$
  for all $\omega\in K$.
  Indeed, we can choose $K$
  with $\LowProbMod(K)$ so close to $1$
  that the sum of $S$
  and a positive capital process eventually attaining
  $2\sup\lvert F\rvert + 2$ on $K^c$
  will give a capital process starting from $\int F\dd W + 2\epsilon$
  and exceeding $F(\omega)$ in the limit.
\end{enumerate}
From now on we fix a compact $K\subseteq\Omega$,
assuming, without loss of generality,
that the statements inside the curly braces
in (\ref{eq:super-modulus}) and (\ref{eq:super-2-variation})
are satisfied for some $\alpha>0$.

In the rest of the proof we will be using,
often following \cite{shafer/vovk:2001}, Section 6.2,
the standard method going back to Lindeberg \cite{lindeberg:1922}.
For $i=N-1$,
define a function
$\overline{U}_i:\bbbr\times[0,\infty)\times\bbbr^i\to\bbbr$
by
\begin{equation}\label{eq:overline-U}
  \overline{U}_i(s,D;s_1,\ldots,s_i)
  :=
  \int_{-\infty}^{\infty} U_{i+1}(s_1,\ldots,s_i,s+z) \Normal_{0,D}(dz),
\end{equation}
where $U_N$ stands for $U$
and $\Normal_{0,D}$ is the Gaussian probability measure on $\bbbr$
with mean $0$ and variance $D\ge0$.
Next define, for $i=N-1$,
\begin{equation}\label{eq:U}
  U_i(s_1,\ldots,s_i)
  :=
  \overline{U}_i(s_i,T/N;s_1,\ldots,s_i).
\end{equation}
Finally, we can alternately use (\ref{eq:overline-U}) and (\ref{eq:U})
for $i=N-2,\ldots,1,0$ to define inductively
other $\overline{U}_i$ and $U_i$
(with (\ref{eq:U}) interpreted as
$
  U_0
  :=
  \overline{U}_0(0,T/N)
$
when $i=0$).
Notice that $U_0=\int F\dd W$.

Informally,
the functions (\ref{eq:overline-U}) and (\ref{eq:U})
constitute Sceptic's goal:
assuming $\omega\in K$,
he will keep his capital at time $iT/N$, $i=0,1,\ldots,N$, close to
$U_i(\omega(T/N),\omega(2T/N),\ldots,\omega(iT/N))$
and his capital at any other time $t\in[0,T]$
close to $\overline{U}_i(\omega(t),D;\omega(T/N),\omega(2T/N),\ldots,\omega(iT/N))$
where $i:=\lfloor Nt/T\rfloor$ and $D:=(i+1)T/N-t$.
This will ensure that his capital at time $T$ is close to $F(\omega)$
when his initial capital is $U_0=\int F\dd W$.

It is easy to check that each function $\overline{U}_i(s,D;s_1,\ldots,s_i)$
satisfies the heat equation in the variables $s$ and $D$:
\begin{equation}\label{eq:heat}
  \frac{\partial\overline{U}_i}{\partial D}
  (s,D;s_1,\ldots,s_i)
  =
  \frac12\frac{\partial^2\overline{U}_i}{\partial s^2}
  (s,D;s_1,\ldots,s_i)
\end{equation}
for all $s\in\bbbr$, all $D>0$,
and all $s_1,\ldots,s_i\in\bbbr$.
This is the key element of the proof.

Sceptic will only bet at the times that are multiples of $T/LN$,
where $L\in\{1,2,\ldots\}$ will later be chosen large.
For $i=0,\ldots,N$ and $j=0,\ldots,L$ let us set
\begin{equation*}
  t_{i,j}
  :=
  iT/N+jT/LN,
  \quad
  S_{i,j}
  :=
  \omega(t_{i,j}),
  \quad
  D_{i,j}
  :=
  T/N-jT/LN.
\end{equation*}
For any array $A_{i,j}$,
we set $dA_{i,j} := A_{i,j+1} - A_{i,j}$.

Using Taylor's formula
and omitting the arguments $\omega(T/N),\ldots,\omega(iT/N)$,
we obtain, for $i=0,\ldots,N-1$ and $j=0,\ldots,L-1$,
\begin{multline}\label{eq:3}
  d \overline{U}_i(S_{i,j},D_{i,j})
  =
  \frac{\partial\overline{U}_i}{\partial s} (S_{i,j},D_{i,j}) dS_{i,j}
  +
  \frac{\partial\overline{U}_i}{\partial D} (S_{i,j},D_{i,j}) dD_{i,j}
\\
  +
  \frac12 \frac{\partial^2\overline{U}_i}{\partial s^2} (S_{i,j}',D_{i,j}')
  (dS_{i,j})^2
  +
  \frac{\partial^2\overline{U}_i}{\partial s \partial D} (S_{i,j}',D_{i,j}')
  dS_{i,j} dD_{i,j}
\\
  +
  \frac12 \frac{\partial^2\overline{U}_i}{\partial D^2} (S_{i,j}',D_{i,j}')
  (dD_{i,j})^2,
\end{multline}
where $(S_{i,j}',D_{i,j}')$ is a point
strictly between $(S_{i,j},D_{i,j})$ and $(S_{i,j+1},D_{i,j+1})$.
Applying Taylor's formula to
$\partial^2\overline{U}_i/\partial s^2$,
we find
\begin{multline*}
  \frac{\partial^2\overline{U}_i}{\partial s^2} (S_{i,j}',D_{i,j}')
  =
  \frac{\partial^2\overline{U}_i}{\partial s^2} (S_{i,j},D_{i,j})\\
  +
  \frac{\partial^3\overline{U}_i}{\partial s^3}
  (S_{i,j}'',D_{i,j}'') \Delta S_{i,j}
  +
  \frac{\partial^3\overline{U_i}}{\partial D \partial s^2}
  (S_{i,j}'',D_{i,j}'') \Delta D_{i,j},
\end{multline*}
where $(S_{i,j}'',D_{i,j}'')$ is a point strictly between
$(S_{i,j},D_{i,j})$ and $(S_{i,j}',D_{i,j}')$,
and $\Delta S_{i,j}$ and $\Delta D_{i,j}$ satisfy
$\lvert\Delta S_{i,j}\rvert \le \lvert dS_{i,j}\rvert$,
$\lvert\Delta D_{i,j}\rvert \le \lvert dD_{i,j}\rvert$.
Plugging this equation and the heat equation (\ref{eq:heat})
into~(\ref{eq:3}), we obtain
\begin{multline}\label{eq:5}
  d\overline{U}_i(S_{i,j},D_{i,j})
  =
  \frac{\partial\overline{U}_i}{\partial s} (S_{i,j},D_{i,j}) dS_{i,j}
  +
  \frac12
  \frac{\partial^2\overline{U}_i}{\partial s^2} (S_{i,j},D_{i,j})
  \left(
    (dS_{i,j})^2 + dD_{i,j}
  \right)
\\
  +
  \frac12
  \frac{\partial^3\overline{U}_i}{\partial s^3}
  (S_{i,j}'',D_{i,j}'')
  \Delta S_{i,j} (dS_{i,j})^2
  +
  \frac12
  \frac{\partial^3\overline{U}_i}{\partial D \partial s^2}
  (S_{i,j}'',D_{i,j}'')
  \Delta D_{i,j} (dS_{i,j})^2
\\
  +
  \frac
    {\partial^2\overline{U}}
    {\partial s \partial D}
  (S_{i,j}',D_{i,j}')
  dS_{i,j} dD_{i,j}
  +
  \frac12
  \frac{\partial^2\overline{U}}{\partial D^2} (S_{i,j}',D_{i,j}')
  (dD_{i,j})^2.
\end{multline}
At this point we are at last ready to define Sceptic's elementary betting strategy:
namely, he plays in such a way
that the increment of his capital between times $t_{i,j}$ and $t_{i,j+1}$ 
is equal to the sum of the first two terms on the right-hand side of (\ref{eq:5}).
Both $\partial\overline{U}_i/\partial s$ and $\partial^2\overline{U}_i/\partial s^2$
are bounded as averages of $\partial U_{i+1}/\partial s$
and $\partial^2 U_{i+1}/\partial s^2$,
and so, eventually,
averages of $\partial U/\partial s$ and $\partial^2 U/\partial s^2$,
respectively.

Let us show that the last four terms
on the right-hand side of (\ref{eq:5})
are negligible when $L$ is sufficiently large
(assuming $T$, $N$, and $U$ fixed).
All the partial derivatives involved in those terms are bounded:
the heat equation implies
\begin{align*}
  \frac{\partial^3\overline{U}_i}{\partial D \partial s^2}
  &=
  \frac{\partial^3\overline{U}_i}{\partial s^2 \partial D}
  =
  \frac12
  \frac{\partial^4\overline{U}_i}{\partial s^4},
\\
  \frac{\partial^2\overline{U}_i}{\partial s \partial D}
  &=
  \frac12
  \frac{\partial^3\overline{U}_i}{\partial s^3},
\\
  \frac{\partial^2\overline{U}_i}{\partial D^2}
  &=
  \frac12
  \frac{\partial^3\overline{U}_i}{\partial D \partial s^2}
  =
  \frac14
  \frac{\partial^4\overline{U}_i}{\partial s^4},
\end{align*}
and $\partial^3\overline{U}_i/\partial s^3$
and $\partial^4\overline{U}_i/\partial s^4$,
being averages of $\partial^3 U_{i+1}/\partial s^3$
and $\partial^4 U_{i+1}/\partial s^4$,
and eventually averages of $\partial^3 U/\partial s^3$
and $\partial^4 U/\partial s^4$,
are bounded.
We can assume that
\begin{equation*}
  \lvert dS_{i,j}\rvert
  \le
  C_1 L^{-1/8},
  \quad
  \sum_{i=0}^{N-1}
  \sum_{j=0}^{L-1}
  (dS_{i,j})^2
  \le
  C_2 L^{1/16}
\end{equation*}
(cf.\ (\ref{eq:super-modulus}) and (\ref{eq:super-2-variation}), respectively)
for $\omega\in K$
and some constants $C_1$ and $C_2$
(remember that $T$, $N$, $U$, and, of course, $\alpha$ are fixed;
without loss of generality we assume that $N$ and $L$ are powers of $2$).
This makes the cumulative contribution of the four terms
have at most the order of magnitude $O(L^{-1/16})$;
therefore,
Sceptic can achieve his goal for $\omega\in K$
by making $L$ sufficiently large.

To ensure that his capital never drops strictly below $\inf F-1$,
Sceptic stops playing as soon as his capital hits $\inf F-1$.
This will never happen when $\omega\in K$
(for $L$ sufficiently large).

\section{Proofs for the L\'evy game}

The following simple lemma will allow us to deduce
the results for the L\'evy game.
\begin{lemma}\label{lem:bounded}
  For each $\alpha>0$ and $T>0$,
  we have
  \begin{equation}\label{eq:bounded}
    \LowProb
    \left\{
      \sup_{t\in[0,T]}
      \lvert\omega(t)\rvert
      \le
      \alpha^{-1/2}
      T^{1/2}
    \right\}
    \ge
    1-\alpha
  \end{equation}
  in the L\'evy game.
\end{lemma}
\begin{proof}
  Starting from initial capital $\alpha$,
  bet $\alpha/T$ on $\omega^2(t)-t$ at time $0$
  and stop playing (set the stake to $0$) at time
  $
    T
    \wedge
    \inf
    \{
      t
      \st
      \lvert\omega(t)\rvert
      =
      \alpha^{-1/2}
      T^{1/2}
    \}
  $;
  the initial capital $\alpha$ will grow to at least
  $\alpha+\frac{\alpha}{T}(\alpha^{-1}T-T)=1$
  if the inner inequality in (\ref{eq:bounded}) is violated.
\end{proof}

It is instructive to compare the bounds
given by the inner inequalities
in (\ref{eq:modulus}) with $\delta:=T$ and in (\ref{eq:bounded}):
the only difference is the constant factor 157 in the former.
Lemma \ref{lem:modulus} itself
will also continue to hold in the L\'evy game;
we will state it in a slightly weakened form:
\begin{lemma}\label{lem:modulus-levy}
  For each $\alpha>0$ and $T>0$,
  \begin{equation*}
    \LowProb
    \left\{
      \forall\delta>0:
      m^T_{\delta}
      \le
      157\,
      \alpha^{-1/2}
      T^{3/8}
      \delta^{1/8}
    \right\}
    \ge
    1-2\alpha.
  \end{equation*}
\end{lemma}
\begin{proof}
  Let
  $
    \tau
    :=
    \inf
    \left\{
      t\ge0
      \st
      \lvert\omega(t)\rvert
      =
      \alpha^{-1/2}
      T^{1/2}
      +
      1
    \right\}
  $
  and let $\omega^{\tau}:[0,\infty)\to\bbbr$ be the stopped $\omega$,
  $\omega^{\tau}(t):=\omega(\tau\wedge t)$.
  Consider the same elementary betting strategies
  as in the proof of Lemma~\ref{lem:modulus}
  except that they now stop playing at time $\tau$.
  Identity (\ref{eq:identity}) with $\omega^{\tau}$ in place of $\omega$
  shows that the corresponding elementary capital processes
  will also be elementary capital processes in the L\'evy game.
  Their combination (analogous to the one in the proof of Lemma~\ref{lem:modulus})
  witnesses that
  \begin{equation*}
    \sup_{t\in[0,T]}
    \lvert\omega(t)\rvert
    \le
    \alpha^{-1/2} T^{1/2}
    \enspace
    \Longrightarrow
    \enspace
    \forall\delta>0:
    m^T_{\delta}
    \le
    157\,
    \alpha^{-1/2}
    T^{3/8}
    \delta^{1/8}
  \end{equation*}
  with lower probability at least $1-\alpha$
  in the L\'evy game;
  it remains to combine this with Lemma~\ref{lem:bounded}.
\end{proof}

In the same way as we obtained
Lemma~\ref{lem:super-modulus} from Lemma~\ref{lem:modulus},
we can now obtain the following corollary from Lemma~\ref{lem:modulus-levy}:
\begin{lemma}\label{lem:super-modulus-levy}
  For each $\alpha>0$,
  \begin{equation*}
    \LowProb
    \left\{
      \forall T\ge1\;
      \forall\delta>0:
      m^T_{\delta}
      \le
      800\,
      \alpha^{-1/2}
      T^{1/2}
      \delta^{1/8}
    \right\}
    \ge
    1-\alpha.
  \end{equation*}
\end{lemma}
\noindent
Lemma~\ref{lem:super-modulus-levy} shows
that Lemma~\ref{lem:tight} will also hold for the L\'evy game.
In a similar way we can get rid of the prime in $\LowProbMod$
in Lemma~\ref{lem:super-2-variation}.

The proof of (\ref{eq:inequality}) for the L\'evy game
proceeds as in the case of the modified L\'evy game:
indeed,
the $7$ reductions in Section \ref{sec:proof}
do not depend on the game being played,
and the elementary betting strategy constructed afterwards
always makes bounded stakes
(see its description after (\ref{eq:5})).

\section{Conclusion}

In this short section we will state two open problems.
First,
what is the class $\AAA$ of all $\UpProb$-measurable subsets of $\Omega$?
It is easy to see that the statement of Theorem \ref{thm:emergence}
that the sets in $\FFF$ are $\UpProb$-measurable
can be strengthened:
\begin{proposition}
  Each set $A\in\FFF^W$ in the completion of $\FFF$ w.r.\ to $W$
  is $\UpProb$-measurable.
\end{proposition}
\begin{proof}
  To establish (\ref{eq:to-show}) for $A\in\FFF^W$
  we choose $A_1,A_2\in\FFF$
  such that $A_1\subseteq A\subseteq A_2$ and $W(A_1)=W(A_2)$,
  and define $F$ as in the proof of Theorem \ref{thm:emergence}
  (see Section~\ref{sec:expectation}).
  Since now
  $\UpExpect(F)<\UpProb(E)+\epsilon$,
  $\III_{E\cap A}\le F\III_{A}\le F\III_{A_2}$,
  and $\III_{E\cap A^c}\le F\III_{A^c}\le F\III_{A_1^c}$,
  it suffices to notice that
  \begin{equation*}
    \UpExpect(F\III_{A_2})
    +
    \UpExpect(F\III_{A_1^c})
    \le
    \UpExpect(F)
  \end{equation*}
  immediately follows from Theorem \ref{thm:expectation}.
\end{proof}
\noindent
In particular, it would be interesting to know
whether $\AAA$ coincides with $\FFF^W$.

The second problem is:
will the (modified) L\'evy game remain coherent
if the measurability restrictions on stopping times and stakes are dropped?
(In other words,
if each $\sigma$-algebra considered
is extended to become closed under arbitrary,
and not just countable, unions and intersections.)
A positive answer would lead
to simpler and more intuitive definitions.
A negative answer would also be of great interest,
providing a counter-intuitive phenomenon
akin to the Banach--Tarski paradox.
A related question is whether dropping the requirement
that $M$ and $V$ should be bounded
will lead to loss of coherence.

\subsection*{Acknowledgments}

This work was partially supported by EPSRC (grant EP/F002998/1),
MRC (grant G0301107),
VLA,
EU FP7 (grant 201381),
and the Cyprus Research Promotion Foundation.

\ifnotJOURNAL

\fi

\end{document}